\documentclass{gtart2}
\usepackage{amssymb, enumerate, epsfig}
\def\<{\langle}
\def\>{\rangle}
\def\|{{\ |\ }}

\newcommand{\inv}{^{-1}}              
\newcommand{\basl}{\backslash}

\newcommand{\C}{{\mathbb C}}
\newcommand{\N}{{\mathbb N}}    

\newcommand{\R}{{\mathbb R}}

\newcommand{\MCGDn}{\mathcal{MCG}(D_n)}
\renewcommand{\phi}{\varphi}

\def\dt{{\, \cdot \,}}
\newtheorem{lemma}{Lemma}[section] 
\newtheorem{theorem}[lemma]{Theorem}
\newtheorem{prop}[lemma]{Proposition}

\newtheorem{lem}[lemma]{Lemma}
\newtheorem{conj}[lemma]{Conjecture}
\newtheorem{remark}[lemma]{Remark}
\newtheorem{example}[lemma]{Example}

\begin{document}
\title{An algorithm for the word problem in braid groups}
\author{Bert Wiest}
\address{UFR de Math\'ematiques, Universit\'e de Rennes 1, Campus de Beaulieu,
\\ 35042 Rennes cedex, France; {\tt bertw@math.univ-rennes1.fr}}
\begin{abstract}
We suggest a new algorithm for finding a canonical representative of a given 
braid, and also for the harder problem of finding a $\sigma_1$-consistent 
representative. We conjecture that the algorithm is quadratic-time. 
We present numerical evidence for this conjecture, and prove two results: 
(1) The algorithm terminates in finite time. (2) The conjecture holds in 
the special case of 3-string braids - in fact, we prove that the algorithm 
finds a minimal-lenght representative for any 3-string braid.
\end{abstract}
\makeshorttitle


\section{Introduction}

In this paper we propose an algorithm for finding a unique 
short representative for any given element of the Artin braid group
$$
B_n \cong \< \sigma_1,\ldots,\sigma_{n-1} \| [\sigma_i,\sigma_j]=0 
\hbox{ if } |i-j|\geqslant 2, \ \sigma_i \sigma_{i+1}\sigma_i =
\sigma_{i+1}\sigma_i \sigma_{i+1} \>.
$$
Several algorithms with the same aim are already
known, for instance Artin's combing of pure braids \cite{Artin}, as
well as Garside's
\cite{Garside,Epstein,FGM} and Dehornoy's \cite{Dehandle} algorithms.
Still, we believe that the new algorithm, which we call the 
\emph{relaxation algorithm,} is of theoretical interest:
if our conjecture that the algorithm is efficient -- a conjecture which
is supported by strong numerical evidence -- is correct, then we would
be dealing with a new type of convexity property of mapping class groups.
A particularly surprising aspect is that the idea of the algorithm 
appears to have been well-known to the experts for decades
(we just have to be specific about some details), only its efficiency 
has apparently been overlooked.

The algorithm is very geometric in nature, exploiting the natural
identification of $B_n$ with $\MCGDn$, the mapping class group of the 
$n$ times punctured disk. The idea is very naive, and not new -- 
c.f.~\cite{Larue,FGRRW}. One thinks of the disk $D_n$ as being made of an
elastic material; if a given braid is represented by a homeomorphism 
$\phi\co D_n \to D_n$, then one obtains a representative of the inverse 
of the braid by authorizing the puncture points to move, and letting the 
map $\phi$ relax into the identity map. However, this relaxation process
is decomposed into a sequence of applications of generators of the
braid group, where in each step one chooses the generator which
reduces the tension of the elastic $D_n$ by as much as possible. Of course,
both the notion of ``tension/relaxation'', as well as the choice of
generating system of the braid group, must be specified carefully.
Here are, in more detail, the properties of the algorithm.

(a) The algorithm appears to be quadratic-time in the length of the
input braid, and the length of the output braid appears to be bounded
linearly by the length of the input braid. Unfortunately, we are not 
currently able to prove this. However, we shall prove some partial 
results, and report on some strong numerical evidence supporting the 
conjecture that the algorithm is efficient in the above sense.

(b) The algorithm can be finetuned to output only $\sigma_1$-consistent
braids (in the sense of Dehornoy \cite{Dehbook}: the output braid words 
may contain the letter $\sigma_1$ but not $\sigma_1\inv$, or vice versa). 
The bounds from (a) on running time and length of the output word still 
seem to hold for this modified algorithm. This is remarkable,
because it is not currently known whether every braid of lenght $l$ admits
a $\sigma_1$-consistent representative of length $c(n)\cdot l$, where
$c(n)>1$ is a constant depending only on the number of strings $n$.

Characteristics (a) and (b) are almost shared by Dehornoy's handle reduction
algorithm \cite{Dehandle}: it is believed to be of cubic complexity
and to yield $\sigma_1$-consistent braids of linearly bounded length, but
this is currently only a conjecture. 

(c) The algorithm generalizes to other surface braid 
groups. We conjecture that in this setting as well it is efficient in
the sense of (a) above (note that we have no reasonable notion of 
$\sigma_1$-consistency here).

(d) The main interest of the algorithm, however, is theoretical: any 
step towards finding a polynomial bound on its algorithmic complexity 
would likely give new insights into the interactions between
the geometry of the Cayley graph of $B_n$ and the dynamical properties
of braids. Indeed, we conjecture that the output braid words are 
quasi-geodesics in the Cayley graph.
Moreover, the problems raised here may be linked to 
the question whether train track splitting sequences are quasi-geodesics,
as well as to the efficiency of the Dehornoy algorithm \cite{Dehandle}.

The plan of this paper is as follows: in section \ref{thealg} we describe
the algorithm and its possible modifications. In section \ref{numerical}
we present numerical evidence that the algorithm is efficient.
In section \ref{terminates} we prove that it terminates in finite time.
In section \ref{3string} we prove that in the special case of 3-string 
braids the algorithm is indeed quadratic time, and finds a shortest 
possible representative for any braid. 


\section{The algorithm}\label{thealg}

\begin{figure}[htb] 
\centerline{
\input{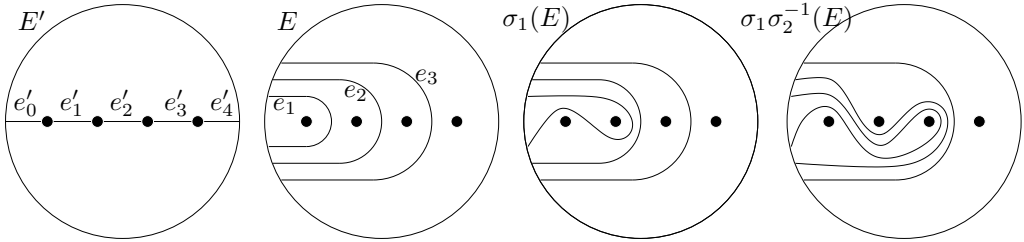} 
}
\caption{The diagram $E'$, and reduced curve diagrams for $id$, $\sigma_1$ 
and $\sigma_1 \sigma_2^{-1}\in B_4$} \label{cddef}
\end{figure}

We start by setting up some notation.
We denote by $D_n$ ($n\in \N$) the closed disk in the complex plane
which intersects $\R$ in the interval $[0,n+1]$, but with the
points $1,\ldots,n \in \C$ removed.
We recall that the braid group $B_n$ is naturally isomorphic to the
mapping class group $\MCGDn$. We denote by $E$ the diagram
in $D_n$ consisting of $n-1$ properly embedded line segments intersecting
the real axis halfway between the punctures, as indicated in figure 
\ref{cddef}. We shall also consider the diagram $E'$ in $D_n$, consisting
of $n+1$ horizontal open line segments. The arcs of both diagrams are 
labelled as indicated.

Given any homeomorphism 
$\phi\co D_n \to D_n$ with $\phi|_{\partial D^2} = id$, 
we write $[\phi]$ for the element of $\MCGDn$ represented by $\phi$.
A \emph{curve diagram} for $[\phi]\in \MCGDn$ is the image $\phi(E)$ of 
$E$ under any homeomorphism $\phi\co D_n \to D_n$ representing $[\phi]$. 
We say a curve diagram is \emph{reduced} if it intersects the diagram $E'$ 
in finitely many points, 
transversely, in such a way that $\phi(E)$ and $E'$ together do not enclose 
any bigons. This is equivalent to requiring that the number of intersection 
points $\phi(E) \cap E'$ be minimal among all homeomorphisms representing 
$[\phi]$ (see e.g.\ \cite{FGRRW}). 
Reduced curve diagrams are essentially unique in the sense
that any two reduced curve diagrams of $[\phi]\in \MCGDn$ can be deformed
into each other by an isotopy of $D_n$ which fixes $E'$ setwise.

We shall call the diagram $E$ the ``trivial curve diagram'', because it
is a reduced curve diagram for the trivial element of $\MCGDn$.

We define the \emph{complexity} of an element $[\phi]$ of $\MCGDn$, and of 
any curve diagram of $[\phi]$, to be the 
number of intersections points of a reduced curve diagram of $[\phi]$
with $E'$. In other words, the complexity of $[\phi]$ is 
$\min |\psi(E) \cap E'|$, where the minimum is taken over all 
homeomorphisms $\psi$ isotopic to $\phi$.
For instance, the elements $[id], \sigma_1$ and 
$\sigma_1 \sigma_2^{-1}$ of $\mathcal{MCG}(D_4)$ have complexity 
$3, 5$ and $9$, respectively. We remark that the complexity
can grow exponentially with the number of crossings of the braid.
(The reader may feel that our definition of complexity is quite arbitrary.
In informal computer experiments we have tested some variations of
this definition, and the results appear to be qualitatively unchanged.)


Next we define the set of generators of the braid group $\MCGDn$ which
will be useful for our purposes: we define a \emph{semicircular move}
to be any element of $\MCGDn$ either of the form 
$\sigma_i^\epsilon \sigma_{i+1}^\epsilon \ldots \sigma_{j-1}^\epsilon
\sigma_j^\epsilon$ with $i<j$, or of the form 
$\sigma_i^\epsilon \sigma_{i-1}^\epsilon \ldots \sigma_{j+1}^\epsilon
\sigma_j^\epsilon$ with $i>j$; in either case we have 
$\epsilon = \pm 1$ and $i,j \in \{1,\ldots, n-1\}$.
In order to explain the name, we remark that semicircular braids can
be realised by a semicircular movement of one puncture of $D_n$ in the 
upper or lower half of $D_n$ back into $D_n \cap \R$. (The movement is
in the upper half if $j>i$ and $\epsilon=1$, or if $j<i$ and $\epsilon=-1$,
it is in the lower half in the reverse case, and when $i=j$ then it can be 
regarded as lying in either half.) There are $2(n-1)^2$ different 
semicircular moves in $\MCGDn$.

We say a semicircular movement is \emph{disjoint from} a curve diagram $D$ 
if it is given by a movement of a puncture along a semicircular arc in the
upper or lower half of $D_n$ which is disjoint from the diagram $D$.

\sh{The ``standard'' relaxation algorithm}
We are now ready to describe the algorithm. One inputs a braid $[\phi]$, 
(e.g.\ as a word in the generators $\sigma_i$), 
and it outputs an expression of $[\phi]^{-1}$ as a product of semicircular
moves (and thus of standard generators $\sigma_i^{\pm 1}$):

{\it Step 1 } Construct the curve diagram $D$ of $[\phi]$. 

{\it Step 2 } For each semicircular movement $\gamma$ which exists in 
$B_n$ and which is disjoint from $D$, calculate the complexity of $\gamma(D)$.

{\it Step 3 } Among the possible moves in step 2, choose the one 
($\gamma'$ say) which yields the minimal complexity. (If several 
moves yield the same minimal complexity, choose one of them arbitrarily or
e.g.\ the lexicographically smallest one.) 
Write down the move $\gamma'$, and define $D := \gamma'(D)$.

{\it Step 4 } If $D$ is the trivial diagram, terminate. If not, go to 
Step 2.

This could be summarized in one sentence: iteratively untangle the curve 
diagram of $[\phi]$, where each iteration consists of the one semicircular 
move which simplifies the curve diagram as much as possible. We shall
prove in section \ref{terminates} that the algorithm terminates, i.e.
that every nontrivial curve diagram admits a move which decreases the 
complexity. Note that there is no
obvious reason to expect this algorithm to be efficient: one might expect
that one has to perform some seemingly inefficient steps first which
then allow a very rapid untangling later on. Yet, as we shall see,
it appears that such phenomena do not occur.

\sh{The $\sigma_1$-consistent version of the relaxation algorithm}
Next we indicate how the algorithm can be improved in order to output
only $\sigma_1$-consistent braid words. We recall that a braid word is 
\emph{$\sigma_1$-consistent} if it contains only the letter $\sigma_1$
but not $\sigma_1^{-1}$ (``$\sigma_1$-positive''), or if it contains 
$\sigma_1^{-1}$ but not $\sigma_1$ (``$\sigma_1$-negative''),
or indeed if it contains no letter $\sigma_1^{\pm 1}$ at all
(``$\sigma_1$-neutral''). 
It is a theorem of Dehornoy \cite{Dehbook} that every braid has a
$\sigma_1$-consistent representative. More precisely, it was shown in
\cite{FGRRW} that a braid $[\phi]$ is $\sigma_1$-positive if and only if in 
a reduced curve diagram $\phi(E)$, the ``first'' intersection of $\phi(e_1)$
with $E'$ lies in $e'_0$, where $e_1$ is oriented from bottom to top.
By contrast, the braid is $\sigma_1$-negative if and only if in the 
opposite orientation of $e_1$ (from top to bottom) the ``first'' intersection
with $E'$ lies in $e'_0$. Finally, a braid is $\sigma_1$-neutral if and only
if $\phi(e_1)\sim e_1$. By abuse of notation we shall also speak of a
reduced curve diagram as being $\sigma_1$-positive, negative, or neutral.
Roughly speaking, a diagram is $\sigma_1$-positive if $\phi(e_1)$ starts
by going up, and $\sigma_1$-negative if it starts by going down.
We remark that these geometric conditions are easy to check once the curve 
diagram of a braid has been calculated: using the notation of the appendix,
a curve diagram is $\sigma_1$-positive if $d^0_/>0$ holds, it is 
$\sigma_1$-negative if $d^0_\basl > 0$ holds, and it is $\sigma_1$-neutral
if $d^0_/ = d^0_\basl = 0$. For instance, one sees in figure~\ref{cddef} that 
the curve diagrams for both
$[\sigma_1]$ and $[\sigma_1 \sigma_2^{-1}]$ are $\sigma_1$-positive.

Here are the changes that need to be made to the above algorithm: after
Step 1 is completed, one checks whether the braid $[\phi]$ is 
$\sigma_1$-positive, negative, or neutral. For definiteness let us say
it is $\sigma_1$-positive (the $\sigma_1$-negative case being symmetric). 
Thus during the untangling process we want to avoid using semicircular 
moves which involve the letter $\sigma_1$, only $\sigma_1^{-1}$s are
allowed. Thus in Step 3  we cannot choose among \emph{all} semicircular 
moves which are disjoint from $D$, but we restrict our choice to those 
moves $\gamma$ which satisfy
\begin{itemize}
\item[(a)] the move $\gamma$ does not involve the letter $\sigma_1$,
\item[(b)] the braid $[\phi]\gamma^{-1}$ is not $\sigma_1$-negative;
in other words, the curve diagram $\gamma(D)$ is still $\sigma_1$-positive
or $\sigma_1$-neutral, but not $\sigma_1$-negative.
\end{itemize}

Note that condition (b) is really necessary: if applying $\gamma$ turned our
$\sigma_1$-positive curve diagram $D$ into a $\sigma_1$-negative one, then
we would have no chance of completing the untangling process without using
the letter $\sigma_1$ later on.

\begin{remark}\rm The restriction to semicircular moves \emph{which are 
disjoint from the curve diagram} is superfluous -- all results in this 
paper are true
without it, and indeed for the standard algorithm the semicircular movement 
which reduces complexity by as much as possible is automatically disjoint
from the diagram. We only insist on this restriction because it simplifies 
the proof in section \ref{3string}
\end{remark}

\begin{example}\rm \label{algex} 
The curve diagram of the braid $\sigma_1\inv\sigma_2\sigma_1$ is shown in 
figure \ref{F:algex}(a). The move $\sigma_1\inv \sigma_2\inv$ reduces the 
complexity of the diagram by 4, whereas $\sigma_2$ (the only other 
semicircular move disjoint from the diagram) reduces the complexity only 
by 2. However, the move $\sigma_1\inv \sigma_2\inv$ is forbidden, since 
its action would turn the positive diagram into a negative one. Thus we 
start by acting by $\sigma_2$, and the resulting diagram is shown in 
figure~\ref{F:algex}(b). In this diagram, the action of 
$\sigma_1\inv\sigma_2\inv$ is legal, and relaxes the diagram into
the trivial one. In total, the curve diagram of 
$\sigma_1\inv\sigma_2\sigma_1$ was untangled by 
$\sigma_2 \; \sigma_1\inv\sigma_2\inv$.
\begin{figure}[htb] 
\centerline{
\input{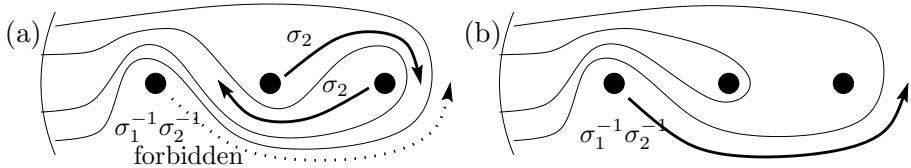} 
}
\caption{Untangling the curve diagram of the braid $\sigma_1\inv\sigma_2\sigma_1$ yields $\sigma_2\dt\sigma_1\inv\sigma_2\inv\dt$}\label{F:algex}
\end{figure}
\end{example}


\section{Numerical evidence}\label{numerical}

The above algorithms were implemented by the author using the programming 
language C \cite{myCprog}. The aim of this section is to 
report on the results of systematic experimentation with this program. 
A brief description how curve diagrams were coded and manipulated can
be found in Appendix 1. Since the algorithm involves calculations with 
integer numbers whose size grows exponentially with the length of the 
input braid, and no software capable of large integer arithmetic was used, 
the search was restricted to braids with up to 50 crossings.

The following notation
will be used: for any braid word $w$ in the letters $\sigma_i^{\pm 1}$ 
we denote by $l(w)$ the number of letters (i.e.\ the number of crossings), 
and by $l_{\rm{out}}(w)$ the number of letters of the output of our algorithm.
The results of the experiments  have two surprising aspects: 
firstly, they strongly support the following

\begin{conj}[Main conjecture]\label{efficcon} 
For every $n>2$ there exists a constant $c(n)\geqslant 1$ such 
that for all braid words $w$ we have $l_{\rm{out}}(w) < c(n) l(w)$. 
\end{conj}

Secondly, they show that the $\sigma_1$-consistent version of the algorithm 
very often yields \emph{shorter} output-braids than the standard version. 

The first table gives the average length of the output braid
among 10.000 randomly generated braid words of length 10, 20, 30, 40
and 50, with 4, 5, and 6 strings. Our random braid words did not
contain any subwords of the form $\sigma_i^\epsilon \sigma_i^{-\epsilon}$,
i.e.\ they did not contain any obvious simplifications. Each entry in
the table is of the form ***/***, where the first number refers to the
standard algorithm, and the second one to the $\sigma_1$-consistent one.

\begin{verbatim} 
       length 10  length 20  length 30  length 40  length 50
  N=4   8.5/8.4   16.0/15.6  23.4/22.8  30.9/29.9  38.3/36.8
  N=5   8.7/8.7   17.2/17.0  25.9/25.3  34.5/33.8  43.2/42.2
  N=6   8.7/8.6   17.3/17.2  26.2/26.2  35.6/35.7  45.0/45.0
\end{verbatim}

The next table shows the worst cases (that is, the longest output words)
that occurred among $\geqslant 40.000.000$ randomly generated braids. 
(Thus every entry in this table required 40.000.000, and in some cases 
much more, complete runs of the algorithm.)

\begin{verbatim} 
       length 10  length 20  length 30  length 40  length 50
  N=4   24 / 20    52 / 38    70 / 54    86 / 66   106 / 80
  N=5   26 / 26    58 / 54    76 / 66   116 / 82   130 / 104
  N=6   30 / 38    78 / 68   102 / 100  140 / 120  178 / 132
\end{verbatim}


While these values appear to support conjecture \ref{efficcon}, one should
keep in mind that the worst examples may be exceedingly rare, and thus
may have been overlooked by the random search. 

If conjecture \ref{efficcon} holds, then we would have that the 
running time of the algorithm (on a Turing machine) is bounded by 
$c'(n) (l(w))^2$, where $c'$ is some other constant depending only on $n$. 
This is because there is an exponential
bound on the numbers $d^i_\supset, d^i_\subset, d^i_{\overline{\cdot\cdot}}, 
d^i_{\underline{\cdot\cdot}}, d^i_/$, and $d^i_\basl$ describing the
curve diagram (see appendix) in terms of the length of the braid word. 
Performing an elementary
operation (addition, comparison etc.) on numbers bounded by
$\exp(l(w))$ takes time at most $O(\log(\exp(l(w))))=O(l(w))$. Since the 
number of operations performed in each step of the algorithm is constant,
and the number of steps which need to be performed is (we conjecture)
bounded linearly by $l(w)$, we get a total running time bounded by
$c'(n) (l(w))^2$.

The algorithm admits variations. For instance, one can use a different
diagram on $D_n$ as a trivial curve diagram (instead of $E$), and one can 
define the
complexity of a curve diagram in a different way (e.g. not counting the
intersection points of $\phi(E)$ but those of $\phi^{-1}(E)$ with the 
horizontal axis). Informal experiments indicate that our results seem
to be ``stable'' under such modifications of the algorithm: computation
time and length of output braids still appear to be quadratic, respectively
linear.

It should also be mentioned that the requirement that the complexity
of the curve diagrams be reduced \emph{as much as possible} in each step
is essential: if one only asks for the complexity to be reduced (by any
amount), then $l_{\rm{out}}(w)$ can depend exponentially on $l(w)$.

The reader should be warned of two attacks on conjecture~\ref{efficcon} 
which appear not to work.
Firstly, running the algorithm on pairs of words $(w,w\sigma_i)$ or 
$(w,\sigma_i w)$, where $\sigma_i$ is a generator, often yields
pairs of output braid words with dramatically different lengths.
Thus it appears unlikely that any fellow-traveller type property between
output braid words can  serve as an explanation for the algorithm's 
efficiency.
Secondly, there are instances of braids where the number of semicircular 
factors in the output braid word is larger than the number of Artin factors
in the input braid word.

One may hope, however, that Mosher's automatic structure on $B_n$ 
\cite{Mosher} can help to prove the main conjecture \ref{efficcon}, 
because the conjecture would follow from a positive answer to the following

{\bf Question } Suppose that the algorithm produces a sequence 
$D^{(1)},\ldots, D^{(N)}$ of curve diagrams of decreasing complexity and
with $D^{(N)}=E$, 
representing braids $[\phi_1],\ldots,[\phi_N]$ in $\mathcal{MCG}(D_n)$. 
(In particular, $[\phi_N]=1$). Is it true that the lengths of the Mosher 
normal forms (the ``flipping sequences'') of these braids are strictly 
decreasing?



\section{The algorithm terminates}\label{terminates}

While we have no explanation for the extraordinarily good performance
of the algorithm, we can at least prove that it terminates in finite time
(and one could easily deduce from our proof an exponential bound on the 
running time). In fact, in order to ascertain the termination of the
algorithm, even of its $\sigma_1$-consistent version, it suffices to prove
\begin{prop}\sl \label{simplexists}
For any reduced, $\sigma_1$-positive curve diagram $D$ 
in $D_n$ there exists a semicircular braid $\gamma$ in the upper half
of $D_n$ such that the diagram $\gamma(D)$ is still $\sigma_1$-positive
or $\sigma_1$-neutral and such that the complexity of $\gamma(D)$ is 
strictly lower than the complexity of $D$.
\end{prop}

We recall that by a semicircular braid ``in the upper half of $D_n$'' we 
mean one which can be realized by a simicircular movement of one of the 
punctures in the part of $D_n$ belonging to the upper half plane.

Let us first see why this proposition implies that the algorithm terminates.
The curve diagram of the input braid $[\phi]$ has finite complexity (and in 
fact if $[\phi]$ is given  as a braid word of length $l$ in the letters 
$\sigma_i^\epsilon$ then one can easily obtain an exponential bound on the 
complexity in terms of $l$ and $n$: the complexity is bounded by
$(n-1)l^3$). If $[\phi]$ is $\sigma_1$-positive, then the proposition 
implies that the curve diagram can be turned into a $\sigma_1$-neutral 
diagram by applying a sequence of semicircular braids which contain no 
letter $\sigma_1$. By a symmetry
argument we can prove an analogue statement for $\sigma_1$-negative
curve diagrams: they can be untangled by applying a $\sigma_1$-positive 
braid. Once we have arrived at a $\sigma_1$-neutral diagram, we can simply
cut $D_n$ along a vertical arc through the leftmost puncture, concentrate
on the part of our cut $D_n$ which contains punctures number $2,\ldots,n$,
and proceed by induction to make the diagram $\sigma_2$-neutral, then 
$\sigma_3$-neutral etc., until it coincides with $E$.

\begin{proof}[Proof of proposition \ref{simplexists}]
Our argument was partly inspired by the unpublished thesis of Larue 
\cite{Larue}. We denote the lower half of $D_n$ by $D_n^\vee$
and the upper half by $D_n^{\wedge}$.
We say a that the $i$th puncture of $D_n$ is \emph{being pushed up} 
by the diagram $D$ if one of the path components of
$D\cap D_n^\vee $ consists of an arc with one endpoint in $e'_{i-1}$ and 
the other in $e'_i$.

We first claim that some puncture other than the first (leftmost) one is
being pushed up by $D$.

To see this, let us denote by $v$ the vertical arc in $D_n$ through the 
first puncture - it is cut into two halves by the first puncture.
Let us write $D_n^>$ for the closure of the component of $D_n\basl v$
which contains punctures number $2,\ldots,n$ (that is, $D_n^>$ equals 
$D_n$ with everything left of $v$ removed). 
Since $D$ is $\sigma_1$-positive, we have that the upper half of $v$ has
two more intersection points with $D$ than the lower half. Therefore
at least one of the components of $D_n^>\cap D$ is an arc both of whose
endpoints lie in the upper half of $v$. Since this arc has to intersect
the diagram $E'$, we deduce that some component of $D\cap D_n^\vee$ has both 
endpoints in intervals $e'_i$, $e'_j$ with $i,j \neq 0$. Now that we have
found at least one ``U-shaped'' arc in $D_n^\vee \cap D_n^>$, we can  
consider an innermost one. Since all punctures of $D_n$ must lie in 
different path components of $D_n\basl D$, all innermost U-shaped arcs 
must have endpoints in a pair of \emph{adjacent} arcs $e'_{i-1}, e'_i$ 
of $E'$ with $i-1 \neq 0$. This completes the proof of the first claim.

It is this $i$th puncture that is going to perform the semicircular
move: our second claim is that there exists an oriented arc $g$ in 
$D_n \basl D$ which starts on the $i$th puncture, lies entirely in 
the upper half of $D_n$, and terminates on one of the horizontal arcs 
$e'_j$ with $j\neq i,i-1$.

To convince ourselves of the second claim, let us look at the two arcs
of $D\cap D_n^\wedge$ which have boundary points in common with our
innermost U-shaped arc in $D\cap D_n^\vee$. Not both of these arcs can
end on $\partial D_n$, so at least one of them, which we shall call $g'$, 
ends on $E'$. Then we can define $g$ to be the arc in $D_n\basl D$ which 
starts on the $i$th puncture, and stays parallel and close to $g'$, until
it terminates on the same arc of $E'$ as $g'$. This construction proves
the second claim.

We now define our semicircular braid $\gamma$ to be the slide of the $i$th
puncture along the arc $g$.

Next, we prove that this move decreases the complexity of the curve diagram. 
Let $D'$ be the diagram obtained from the diagram $D$ by first applying the 
homeomrphism $\gamma$ and then reducing the resulting diagram with respect 
to the horizontal line $E'$.  In comparison to $D$, in this new diagram $D'$
the arc $e'_j$ has been divided into two arcs by the formerly $i$th puncture, 
whereas the arcs $e'_{i-1}$ and $e'_i$ have merged. Moreover, if the diagram 
$D$ contained $k$ U-shaped arcs in $D_n^{\vee}$ with one endpoint
in $e'_{i-1}$ and the other in $e'_i$, then the complexity of $D'$ is
$2k$ lower than the complexity of $D$: our semicircular move $\gamma$
eliminated these $2k$ intersections of $D$ with $E'$, and did not create
or eliminate any others. 

The only thing left to be seen is that the diagram $D'$ is still 
$\sigma_1$-positive or $\sigma_1$-neutral, but not $\sigma_1$-negative.
But a $\sigma_1$-negative diagram could only be obtained 
if the arc $g$ ended in the
interval $e'_0$ (i.e.\ if the last letter of $\gamma$ was $\sigma_1^{-1}$).
The proof is now exactly analogue to the proof of Proposition 4.3, Claim 2 
in \cite{FGRRW} which states that curve diagram obtained by ``sliding a 
puncture along a useful arc'' is not $\sigma_1$-negative.
\end{proof}

%
%


\section{The special case of 3-string braids}\label{3string}

The reader may have wondered why no experimental results for 3-string
braids were presented in section \ref{numerical}. The reason is that we 
can \emph{prove} that the best possible results hold in this case:
\begin{theorem}\sl\label{3stringthm} 
For three-string braids, the $\sigma_1$-consistent version of our algorithm
outputs only braid words of minimal length.
\end{theorem}

Thus for $n=3$, conjecture~\ref{efficcon} holds with $c(n)=1$.
In other words, if we input any braid word $w'$ in the letters 
$\sigma_i^{\epsilon}$ with $i\in \{1,2\}$ and $\epsilon=\pm 1$, then
the output word $w$ satisfies $l(w)\leqslant l(w')$. Unfortunately, the
proof below is very complicated, and it would be nice to have a more
elegant proof.
We remark that the analogue of theorem \ref{3stringthm} for the standard 
algorithm also holds. We shall not prove this here, because the argument 
becomes even more complicated. 

In what follows, we shall always denote by $w$ a braid word output by
the algorithm.
The word $w$ is given as a product of semicircular braids, 
and we shall indicate this decomposition, when necessary, by separating
the semicircular factors by a dot ``$\dt$''. Thus 
$\sigma_2 \sigma_1 \dt \sigma_2^{-1}$ denotes the product of the two 
semicircular braids $\sigma_2 \sigma_1$ and $\sigma_2^{-1}$, and the word
$\sigma_2 \sigma_1^{-1}\dt \sigma_2$ (without a dot between $\sigma_2$ and 
$\sigma_1^{-1}$) is meaningless. Finally, we shall denote by an asterisk
$*$ any symbol among $\{1,-1,2,-2,\dt \}$, or the absence of symbols.

\begin{lem}\sl \label{3strkeylem}
Let $w$ be the output braid word of the $\sigma_1$-consistent algorithm. Then
\begin{itemize}
\item[(i)] $w$ does not contain any subword of the form
$ *\,\sigma_i\dt \sigma_i\inv *$ with $i\in \{1,2\}$.
\item[(ii)] The only place in $w$ where the subword 
$\dt \sigma_2^{-1} \dt \sigma_1 \sigma_2$ can occur is near the end of $w$, 
in the context of a terminal subword $\dt \sigma_2^{-1} \dt \sigma_1 
\sigma_2 (\dt \sigma_2)^k$ with $k\in \N \cup \{0\}$.
\item[(iii)] $w$ does not contain any subword of the form
$\dt\sigma_1\dt\sigma_2\dt$ or $\dt\sigma_2\dt\sigma_1\dt$.
\end{itemize}
Items (i),(ii), and (iii) are also true if every single letter in their 
statement is replaced by its inverse.
\begin{itemize}
\item[(iv)] Subwords of the form $\dt \sigma_2 \sigma_1\dt\sigma_2 \sigma_1\dt
\sigma_2\inv *$ cannot occur in $w$.
\item[(v)] Subwords of the form $\dt\sigma_2 \sigma_1\dt\sigma_2 \sigma_1\dt
\sigma_1*$ cannot occur in $w$.
\item[(vi)] Subwords of the form $\dt \sigma_1\dt \sigma_2\sigma_1\dt$ and
$\dt \sigma_2\sigma_1\dt \sigma_2\dt$ can occur in $w$, but only in the 
context of a \emph{terminal} subword which is positive (i.e.\ has
not a single letter $\sigma_1\inv$ or $\sigma_2\inv$.) Moreover, they cannot
be immediately preceded by the subword $\dt\sigma_2^{-1}$.
\end{itemize}
Items (iv)--(vi) are also true if every single letter in their statement
is replaced by its inverse or if the r\^oles of $\sigma_1$ and $\sigma_2$
are interchanged.
\end{lem}
Items (iii)--(vi) state essentially that the subwords $*\sigma_1\dt\sigma_2*$ 
and $*\sigma_2\dt\sigma_1*$ cannot occur, except in a few very special cases 
near the end of the word.

\begin{proof}[Proof of lemma \ref{3strkeylem}]
For the proof of (i) we suppose, for a contradiction, that $w$ contains a 
subword consisting of two semicircular braids: the first subword consists of
or terminates in $\sigma_i^\epsilon$, the second one consists of or begins 
with $\sigma_i^{-\epsilon}$.
For definiteness we suppose that we have the subword $*\sigma_2\dt 
\sigma_2\inv*$, the other cases being exactly analogue. The first factor
$\sigma_1\sigma_2$ or $\sigma_2$ can be realized by a semicircular movement
in the upper half of $D_n$. We consider the curve diagram $D$ that is
obtained after the action of the first factor. In this diagram, the third
puncture is either being pulled down, or it may be completely to the right
of $D$ (in this case, $D$ can be untangled by a braid which contains only 
the letters $\sigma_1^{\pm 1}$), but it is certainly not being pushed up.
But in this situation, acting on $D$ by $\sigma_2\inv$ or 
$\sigma_2\inv \sigma_1\inv$ cannot reduce the complexity.
(It is worth remarking that the statement (i) is wrong for braids with more
than three strings -- there the first factor might be followed by a 
semicircular movement of the second puncture in the lower half of $D_n$ to
a point further to the right.)


\begin{figure}[htb]  
\centerline{
\input{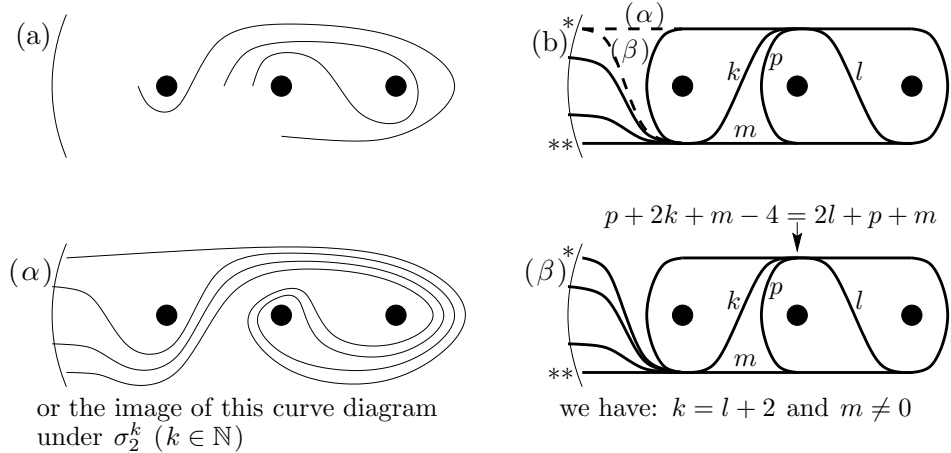} 
}
\caption{Which curve diagrams might give rise to the subword $\sigma_2\inv\dt
\sigma_1\sigma_2$?}
\label{F:trtr}
\end{figure}

The proof of (ii) is more complicated. 
Suppose that $D$ is a curve diagram which the $\sigma_1$-consistent algorithm 
wants to untangle by first performing a $\sigma_2\inv$, and then a 
$\sigma_1 \sigma_2$. Then the curve diagram $D$ must be negative (because 
it can be untangled by a $\sigma_1$-positive braid word), and it must
contain arcs as indicated
in figure \ref{F:trtr}(a) (note that there may be a number of parallel
copies of these arcs). We deduce that $D$ is carried by one of two
possible train tracks, which are both indicated in figure \ref{F:trtr}(b)
-- the two possibilities are labelled $(\alpha)$ and $(\beta)$. The 
labels $k,l,m,p$ are, of course, integer variables which indicate how often 
the corresponding piece of track is traversed by the curve diagram.

Case $(\alpha)$ is quite easy to eliminate: the essential observation is
that in case $(\alpha)$ we must have $m=0$.
To see this, consider the arc of the curve diagram whose extremities are 
the points labelled $*$ and $**$. This arc, together with the segment of
$\partial D_n$ shown, must bound a disk which contains exactly two of the 
punctures. If we had $m\neq 0$, however, then the arc from $*$ to $**$ would 
intersect $E'$ only once, namely in $e'_3$, and the disk would contain all 
three punctures, which is absurd. Thus, in case $(\alpha)$ we must have $m=0$.

It is now an exercice to show that the only curve diagrams carried 
by the train track of figure \ref{F:trtr}(b)($\alpha$) with $m=0$ and 
containing arcs as indicated in figure
\ref{F:trtr}(a) belong to the family of curve diagrams given in figure 
\ref{F:trtr}($\alpha$). However, the algorithm chooses to untangle these 
diagrams with the movement $(\sigma_2\inv\dt)^{2+k}\sigma_1$, not with 
$\sigma_2\inv\dt\sigma_1 \sigma_2...$, as hypothesised. Thus the curve 
diagram $D$ cannot be carried by the train track from figure 
\ref{F:trtr}(b)($\alpha$).

For case $(\beta)$ we consider again the arc of the curve diagram with
endpoints $*$ and $**$. This arc, together with the segment of
$\partial D_n$ shown, must bound a disk which contains two of the punctures.
For this to happen, however, we need that starting from the point $**$ and
following the arc of the diagram we traverse the track labelled $m$, not 
the track labelled $k$. Therefore in case $(\beta)$ we must have $m\neq 0$.

Moreover, the calculation shown in the figure proves that the curve diagram 
$D$ traverses the track labelled $k$ twice more than the track labelled $l$.
This implies that the third puncture is being pushed up by $l$ U-shaped arcs 
below it
(and similarly there are $l$ caps above the second puncture), 
but there are at least $l+1$ U-shaped arcs below the first puncture. 

Thus the standard algorithm would choose to perform the movement 
$\sigma_1 \sigma_2$, and the $\sigma_1$-consistent algorithm only avoids
doing so because the movement $\sigma_1 \sigma_2$ is illegal (it would 
render the diagram positive). However, we know by hypothesis that after 
applying $\sigma_2\inv$, the move $\sigma_1 \sigma_2$ \emph{is} legal.

\begin{figure}[htb]  
\centerline{
\input{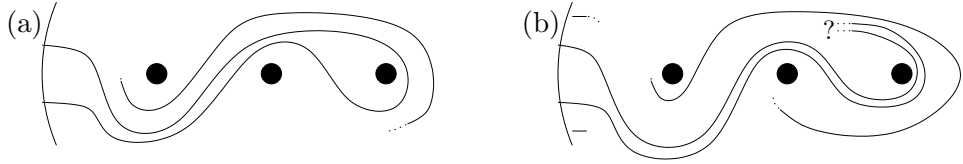} 
}
\caption{Applying $\sigma_2\inv$ to these diagrams legalizes the move 
$\sigma_1\sigma_2$}
\label{F:end3braid}
\end{figure}

This means that the curve diagram is of one of two shapes shown in figure
\ref{F:end3braid}. We observe that the diagram in figure \ref{F:end3braid}(b)
cannot be extended to a reduced curve diagram -- indeed, the reader will find
that any attempt to extend the arcs labelled by a question-mark will yield
arcs of infinite length. Therefore we must be in case (a).
It follows that after applying the braid $\sigma_2\inv\dt\sigma_1\sigma_2$
to the diagram $D$ we are left with a curve diagram that can be untangled 
by $\sigma_2^k$ for some $k\in \N\cup \{0\}$. This completes the proof of (ii).

For the proof of (iii) we remark that the algorithm prefers the
word $\dt \sigma_1^\epsilon \sigma_2^\epsilon\dt$ to the word 
$\dt\sigma_1^\epsilon\dt\sigma_2^\epsilon\dt$, and similarly the word
$\dt\sigma_2^\epsilon \sigma_1^\epsilon\dt$ to the word 
$\dt\sigma_2^\epsilon\dt\sigma_1^\epsilon\dt$ for $\epsilon = \pm 1$.

For the proof of (iv) and (v) we consider a diagram $D$ such that the 
untangling of $D$ begins, according to our algorithm, with the movements
$\sigma_2 \sigma_1\dt \sigma_2 \sigma_1$. We observe that the diagram $D$
must contain curve segements as shown in figure \ref{F:casede}(a) (solid 
lines). The action of $\sigma_2 \sigma_1\dt \sigma_2 \sigma_1$ on
this diagram is indicated by the arrows.

\begin{figure}[htb]  
\centerline{
\input{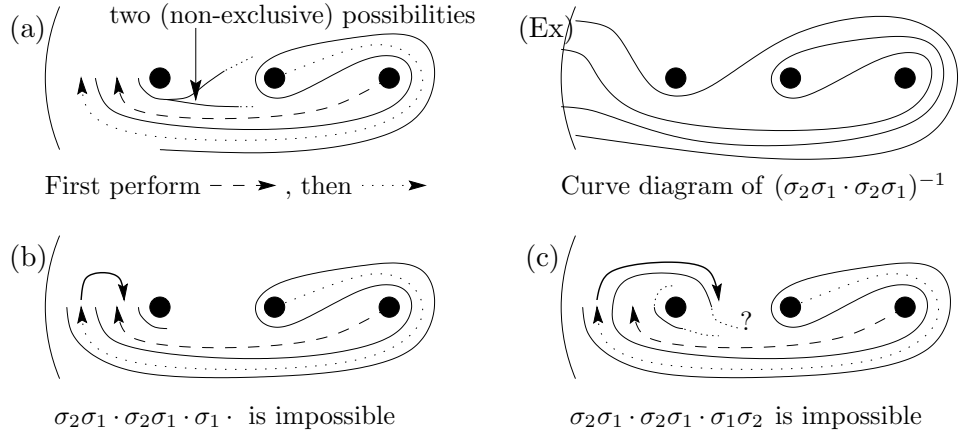} 
}
\caption{Which curve diagrams are relaxed by $\sigma_2 \sigma_1\dt\sigma_2 \sigma_1\dt *$?}
\label{F:casede}
\end{figure}

Let us first prove (iv), i.e.\ that the movement
$\sigma_2 \sigma_1\dt \sigma_2 \sigma_1\dt \sigma_2\inv *$ cannot occur.
After the action of $\sigma_2 \sigma_1\dt \sigma_2 \sigma_1$ on $D$, the 
puncture which started out in the leftmost position is now the rightmost one. 
In order to get the move $\sigma_2\inv$ or $\sigma_2\inv \sigma_1\inv$ next, 
we would need that in the diagram 
$\sigma_2 \sigma_1\dt \sigma_2 \sigma_1(D)$ the rightmost puncture is being
pushed up by a U-shaped arc below it.
However, one can see
that applying $\sigma_2 \sigma_1\dt \sigma_2 \sigma_1$ to the
diagram in figure \ref{F:casede}(a) can never yield such a diagram.

For the proof of (v), we distinguish two cases. First, the sequence of moves
$\sigma_2 \sigma_1\dt \sigma_2 \sigma_1\dt \sigma_1\dt$ is impossible,
because the fat arrow in figure \ref{F:casede}(b) must intersect the
curve diagram $D$. Regarding the second case, if we suppose that the 
curve diagram $D$ gives rise to the moves 
$\sigma_2 \sigma_1\dt \sigma_2 \sigma_1\dt \sigma_1\sigma_2$, then we can add
information to figure  \ref{F:casede}(a): in this case, the diagram $D$ must
have contained the arcs shown in figure \ref{F:casede}(c). However, this
diagram cannot be extended to a reduced curve diagram which is disjoint from 
the dashed arrow: indeed, any attempt to extend the arc labelled by a
question-mark yields an arc of infinite length. This completes the proof 
of (v).

Finally, we turn to the proof of (vi). Again, we suppose that the
algorithm chooses to untangle the curve diagram $D$ by a sequence of
semicircular moves that starts with $\dt \sigma_1\dt \sigma_2\sigma_1\dt$ 
or with $\dt \sigma_2\sigma_1\dt \sigma_2\dt$. We observe that in both cases
the curve diagram must be carried by one of the two train tracks ($\alpha$
or $\beta$) both shown in figure \ref{F:casef}(a), where the tracks labelled
$k$ and $l$ must each be traversed by at least one curve segment of $D$. 
However, as explained in the figure, case ($\alpha$) is completely 
impossible, and case ($\beta$) is possible only if each of the two pieces
of track labelled $k$ and $l$ is traversed by exactly one segment of $D$.

\begin{figure}[htb]  
\centerline{
\input{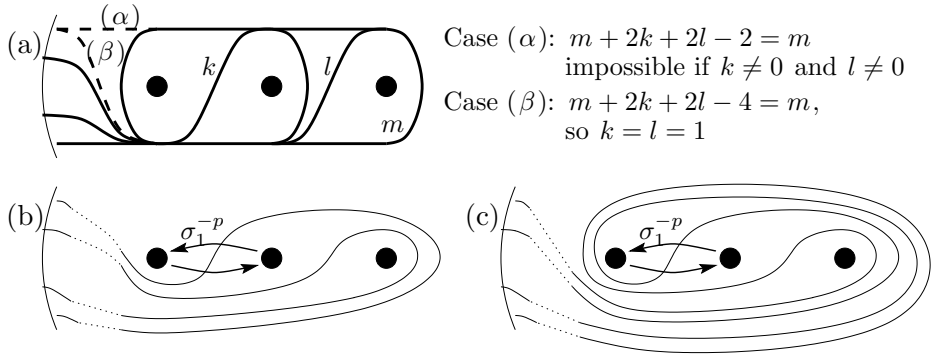} 
}
\caption{Case (vi) of Lemma \ref{3strkeylem}}
\label{F:casef}
\end{figure}

We leave it as an exercice to show that the only curve diagrams fitting these
restrictions are the two shown in figure \ref{F:casef}(b) and (c), 
and all their images under the action of $\sigma_1^{-p} \Delta^{-2q}$, where
$p,q\in \N\cup \{0\}$. In other words, our curve diagram must in fact be
the curve diagram of a braid $\sigma_1\inv \sigma_2\inv \sigma_1\inv 
\sigma_1^{-p} \Delta^{-2q}$ or $\sigma_2\inv \sigma_1\inv \sigma_2\inv 
\sigma_1\inv \sigma_2\inv \sigma_1\inv \sigma_1^{-p} \Delta^{-2q}$, where
$p,q\in \N\cup\{0\}$. Now we observe that any untangling of any of these
curve diagrams that our algorithm may find has only positive letters.
Moreover, if we act by $\sigma_2$ on any of these curve diagrams, then the
relaxation (according to our algorithm) of the resulting diagram starts
with $\sigma_1 \sigma_2$, \emph{not} with $\sigma_2\inv$. This completes
the proof of (vi).
\end{proof}


\begin{proof}[Proof of theorem \ref{3stringthm}]
Let us denote by $\widetilde{w}$ the braid word obtained from an output 
braid word $w$ by removing all the dot-symbols.
Then lemma \ref{3strkeylem} implies that $\widetilde{w}$ cannot 
contain any subword of the form 
$\sigma_2 (\sigma_1 \sigma_2 \sigma_2 \sigma_1)^k \sigma_2^{-1}$, or
$\sigma_2^{-1} (\sigma_1 \sigma_2 \sigma_2 \sigma_1)^k \sigma_2$, or
$\sigma_1 (\sigma_2 \sigma_1 \sigma_1 \sigma_2)^k 
\sigma_2 \sigma_1 \sigma_2^{-1}$, or 
$\sigma_2^{-1} (\sigma_1 \sigma_2 \sigma_2 \sigma_1)^k
\sigma_1 \sigma_2 \sigma_1$, $(k\in \N\cup \{0\})$,
or any of the images of such braids under one of the automorphisms of $B_3$
$(\sigma_1 \to \sigma_2, \sigma_2 \to \sigma_1)$, or 
$(\sigma_1 \to \sigma_1^{-1},  \sigma_2 \to \sigma_2^{-1})$, or
$(\sigma_1 \to \sigma_2^{-1}, \sigma_2 \to \sigma_1^{-1})$.

That is, the word $\widetilde{w}$ admits no obvious simplifications. 
It is well-known that 3-string braids without this type of obvious
shortenings are in fact of minimal length. 
This completes the proof of theorem \ref{3stringthm}.
\end{proof}

\vspace{.4cm}

\begin{appendix}
{\Large{\bf Appendix: Coding of curve diagrams}}

{\small \vspace{0.2cm}
If we cut the disk $D_n$ along $n$ vertical lines through the puncture
points, we obtain $n+1$ vertical ``bands''. The connected components of 
the intersection of a curve diagram with every such band come in six
different types, which are in an obvious way described by the symbols
$\supset, \subset, \overline{\cdot\cdot}, \underline{\cdot\cdot}, /$
and $\basl$.
Thus curve diagrams in $D_n$ are coded by $6(n+1)$ integer numbers 
$d^i_\supset, d^i_\subset, d^i_{\overline{\cdot\cdot}}, 
d^i_{\underline{\cdot\cdot}}, d^i_/$, and $d^i_\basl$ (with $i=0,\ldots n$).
For instance, in figure \ref{F:algex}(b), we have $d^0_/=2, 
d^0_{\overline{\cdot\cdot}}=2, d^1_\basl=2,d^1_{\overline{\cdot\cdot}}=2,
d^2_{\overline{\cdot\cdot}}=1, d^2_\supset=1, d^2_{\underline{\cdot\cdot}}=1,
d^3_\supset=1$, and all other coeffiecients equal zero.

Using this coding, reduced curve diagrams of arbitrary braids can be
calculated by the rule that applying $\sigma_i$ to a curve diagram $D$
changes the diagram as follows. Firstly, all $d_\times^j$ (where $\times$ 
can be any symbol) with $j\notin \{i-1, i, i+1\}$ are unchanged. Apart from
that, the following rules apply simultaneously:

$d^{i-1}_\supset \leftarrow d^{i-1}_\supset$\\
$d^{i-1}_\subset \leftarrow d^i_\subset + d^i_{\basl} + \max(0,d^i_{\underline{\cdot\cdot}}-d^{i-1}_{\underline{\cdot\cdot}}-d^{i-1}_{\basl})$\\
$d^{i-1}_{\overline{\cdot\cdot}} \leftarrow d^{i-1}_{\overline{\cdot\cdot}} + d^{i-1}_\basl - \min(d^{i-1}_\subset, \max(0, d^i_{\underline{\cdot\cdot}}-d^{i-1}_{\underline{\cdot\cdot}}))$\\
$d^{i-1}_{\underline{\cdot\cdot}} \leftarrow \min(d^{i-1}_{\underline{\cdot\cdot}}, d^i_{\underline{\cdot\cdot}})$ \\ 
$d^{i-1}_/ \leftarrow d^{i-1}_/ + d^{i-1}_{\underline{\cdot\cdot}} -
\min(d^{i-1}_{\underline{\cdot\cdot}}, d^i_{\underline{\cdot\cdot}})$\\
$d^{i-1}_\basl \leftarrow \min(d^{i-1}_\subset, \max(0, d^i_{\underline{\cdot\cdot}}-d^{i-1}_{\underline{\cdot\cdot}}))$\\
$d^i_\supset \leftarrow \min(\min(d^i_/, d^{i+1}_\supset) , \max(0,d^{i-1}_\subset + d^{i-1}_{\underline{\cdot\cdot}} - d^i_{\underline{\cdot\cdot}}))$ \\
$d^i_\subset \leftarrow \min(\min(d^{i-1}_\subset, d^i_/),\max(0,d^{i+1}_{\overline{\cdot\cdot}} + d^{i+1}_\basl - d^i_{\overline{\cdot\cdot}}))$ \\
$d^i_{\overline{\cdot\cdot}} \leftarrow d^i_\supset d^i_\basl + d^i_{\overline{\cdot\cdot}} - \min(d^{i-1}_\subset, d^i_/, d^i_{\underline{\cdot\cdot}})$ \\
$d^i_{\underline{\cdot\cdot}} \leftarrow d^i_\subset + d^i_\basl + d^i_{\underline{\cdot\cdot}} - \min(d^{i+1}_\supset , d^i_/, d^i_{\underline{\cdot\cdot}})$ \\
$d^i_/ \leftarrow \max(0,\min(d^i_/, d^{i-1}_\subset)+\min(d^i_/,d^{i+1}_\supset) - d^i_/ )$ \\
$d^i_\basl \leftarrow d^i_\supset+d^i_\subset+d^i_\basl+\max(0,d^i_/ - d^{i-1}_\subset-d^{i+1}_\supset)$ \\
$d^{i+1} \leftarrow d^i_\supset+d^i_\basl+\max(0,d^i_{\overline{\cdot\cdot}}-d^{i+1}_{\overline{\cdot\cdot}}-d^{i+1}_\basl$ \\
$d^{i+1}_\subset \leftarrow d^{i+1}_\subset$ \\
$d^{i+1}_{\overline{\cdot\cdot}} \leftarrow \min(d^{i+1}_{\overline{\cdot\cdot}},d^i_{\overline{\cdot\cdot}})$ \\
$d^{i+1}_{\underline{\cdot\cdot}} \leftarrow d^{i+1}_{\underline{\cdot\cdot}} + d^{i+1}_\basl - \min(d^{i+1}_\basl,\max(0,d^i_{\overline{\cdot\cdot}}-d^{i+1}_{\overline{\cdot\cdot}}))$ \\
$d^{i+1}_/ \leftarrow d^{i+1}_/ + d^{i+1}_{\overline{\cdot\cdot}} - \min(d^{i+1}_{\overline{\cdot\cdot}},d^i_{\overline{\cdot\cdot}})$ \\
$d^{i+1}_\basl \leftarrow \min(d^{i+1}_\basl,\max(0,d^i_{\overline{\cdot\cdot}}-d^{i+1}_{\overline{\cdot\cdot}}))$ \\
} 
\end{appendix}\\
\vspace{2mm}

{\bf Acknowledgement } I thank Lee Mosher and Juan Gonzalez-Meneses for
helpful discussions, and Heather Jenkins and Sandy Rutherford of the
Pacific institute for the mathematical sciences (PIMS) at UBC Vancouver
for letting me use their most powerful computers.

\end{document}